# A FOURIER TRANSFORM METHOD FOR NONPARAMETRIC ESTIMATION OF MULTIVARIATE VOLATILITY

By Paul Malliavin and Maria Elvira Mancino[1]

*Académie des Sciences, Institut de France and University of Firenze*

We provide a nonparametric method for the computation of instantaneous multivariate volatility for continuous semi-martingales, which is based on Fourier analysis. The co-volatility is reconstructed as a stochastic function of time by establishing a connection between the Fourier transform of the prices process and the Fourier transform of the co-volatility process. A nonparametric estimator is derived given a discrete unevenly spaced and asynchronously sampled observations of the asset price processes. The asymptotic properties of the random estimator are studied: namely, consistency in probability uniformly in time and convergence in law to a mixture of Gaussian distributions.

**1. Introduction.** The volatility is a key parameter in financial economics-mathematics and the recent financial econometrics literature has devoted much attention to its computation. As a matter of fact, volatility and co-volatility computation is crucial in hedging strategies (in the classical Black–Scholes environment or in more sophisticated stochastic volatility models), in risk management (e.g., VaR methodologies), in forecasting (e.g., ARCH-GARCH models) and for optimal portfolio selection.

Volatility can be computed through parametric or nonparametric methods, see, for instance, the review by [4]. In the first case, the expected volatility is modeled through a functional form of market or latent variables. Nonparametric methods address the computation of the historical volatility without assuming a functional form of the volatility. As volatility changes over time, its computation through nonparametric methods concentrates on a small time window (a day, a week) and high frequency data are employed.

Received November 2006; revised January 2008.
[1]Supported by the MIUR Grant 206132713-001.
*AMS 2000 subject classifications.* Primary 62G05, 62F12, 42A38; secondary 60H10, 62P20.
*Key words and phrases.* Continuous semi-martingale, instantaneous co-volatility, nonparametric estimation, Fourier transform, high frequency data.

In this paper we propose a nonparametric estimation methodology based on Fourier analysis which applies to continuous semi-martingales and is mainly designed for measuring instantaneous multivariate volatility, exploiting high frequency observations.

Suppose the time evolution of the asset price is *a Brownian semi-martingale* of the form

$$dp(t) = \sigma(t, W) \, dW(t) + b(t, W) \, dt,$$

where $W$ is a Brownian motion. Using Itô calculus, the instantaneous volatility of the process $p$ is obtained

$$(1) \qquad \text{Vol}(p)(t) := \sigma^2(t) = \lim_{\varepsilon \to 0+} E^{\mathcal{F}_t}\left(\frac{(p(t+\varepsilon) - p(t))^2}{\varepsilon}\right).$$

Nevertheless it is difficult to use this formula for computing numerically the instantaneous volatility. As a matter of fact volatility measurement takes place over discrete time intervals and only a single path of the market price is effectively observed; therefore, the conditional expectation $E^{\mathcal{F}_t}$, with respect to the past information $\mathcal{F}_t$, cannot be computed from the observations. As a consequence, volatility is mainly computed over discrete time intervals (*integrated volatility*), relying upon the *quadratic variation formula*. In fact a classical result, essentially due to Wiener, states the following formula which holds almost surely

$$\text{QV}|_{t_0}^{t_1}(p) = \int_{t_0}^{t_1} \sigma^2(s) \, ds,$$

where

$$\text{QV}|_{t_0}^{t_1}(p) := \lim_{n \to \infty} \sum_{0 \leq k < (t_1 - t_0) 2^n} (p(t_0 + (k+1)2^{-n}) - p(t_0 + k 2^{-n}))^2.$$

In the recent financial-econometrics literature many results exploit the quadratic variation formula obtaining nonparametric estimators of integrated volatility and co-volatilities. The *realized volatility—quadratic variation* estimators have been intensively studied and used for financial-econometrics purposes in a series of papers, among them [5, 8, 25, 43]. In order to handle irregularly, even randomly spaced data, the quadratic variation theorem was extended [2, 9, 28, 37].

Nevertheless when instantaneous volatility and co-volatilities are computed by exploiting high frequency data, according to quadratic variation and co-variation formulae, three bottlenecks appear.

The *first* one is specific to the multivariate setting: returns are recorded at the highest available observation frequency; therefore, they are asynchronous across different assets. The realized covariance type estimators require the choice of a "synchronization" method, thus these estimators suffer from a downward bias, when the sampling interval is reduced (known as *Epps effect*, by [16]). Reference [42] provides an analytical study of the realized co-volatility estimator in a general framework which includes asynchronous trading; recently [20] have proposed an estimator of integrated co-volatility which is consistent under asynchronous observations.

*Second*, instantaneous volatility computation (1) involves a sort of numerical derivative. Therefore, the authors, who address the problem of estimating the instantaneous volatility, use a double asymptotic in order to perform both the numerical derivative involved in formula (1) and the approximation procedure (see [14, 18, 37]). However, the empirical derivative required to get the instantaneous volatility gives rise to strong numerical instabilities; moreover, none of these authors address the problem of asynchronous observations in the estimation of instantaneous co-volatility. A different approach through wavelets' series is developed in [19] for deterministic volatility function. Finally, functional methods, which are local in space, have been developed for estimating the volatility as function of the underlying state variable level; see [1, 7, 17, 26].

*Third*, the limiting argument $n \to \infty$ relying on high frequency data cannot be effectively accurate due to market microstructure effects (such as discreteness of prices, bid/ask bounce, etc.), which cause the discrepancy between asset pricing theory based on semi-martingales and the data at very fine intervals. The impact of market microstructure noise on the realized volatility estimator has been studied, among others, in [3, 15, 41, 44], and modifications of realized variance estimator have been proposed to correct the bias due to microstructure noise, for example, the two-scales estimator by [43] and the realized kernels by [10].

The estimation procedure proposed in this paper is immune from the first and second bottlenecks highlighted above, in virtue of its own definition. Fourier estimator uses all the available observations and avoids any "synchronization" of the original data, because it is based on the integration of the time series of returns rather than on its differentiation. As well as regarding the third bottleneck, which concerns the robustness of the estimation procedure under microstructure noise, the Fourier estimator is competitive, in terms of mean squared error, with the methods specifically designed to handle market microstructure noise with only a slightly higher bias; see [36].

The first step in the Fourier analysis methodology for the computation of co-volatility is [32], where the instantaneous multivariate volatility function is expanded into trigonometric polynomials, whose coefficients depend on the log-return processes. The efficiency of the Fourier series procedure for computing and forecasting the *integrated* volatility has been analyzed in [12, 13, 21, 29, 38]. Moreover, an application of Fourier method to the computation of correlation dynamics of futures markets is performed in [24].

In this paper we extend the results of [32]: we prove a general identity which relates the Fourier transform of the co-volatility function with the Fourier transform of the log-returns under the hypothesis that the volatility process is square integrable. From this identity we derive an estimator of *instantaneous* co-volatility based on a discrete, unevenly spaced and asynchronously sampled asset prices. Finally, we prove the consistency in probability uniformly in time of the proposed estimator and the convergence in law to a mixture of Gaussian distribution. As a byproduct, the *integrated* volatility and co-volatilities estimators are obtained by considering the $k = 0$ coefficient of the Fourier expansion. Our theory implies that the Fourier estimator of integrated co-volatility is consistent under asynchronous observations in the absence of microstructure noise.

It is worth noting that the Fourier estimator has two adjustment parameters: the number of coefficients computed for the price and the number of Fourier coefficients used for the reconstruction of the volatility. Balancing the choice of these two parameters according specific situations, the estimator can act as an optical microscope with different magnifications: in the case where microstructure is supposed to be produced by irrelevant noise, low magnification will be used. In both cases, the whole set of observed data are employed. In summary, the Fourier estimator is designed specifically for high frequency data: by cutting the highest frequencies, it uses as much as possible of the sample path without being more sensitive to market frictions; see also [36, 38].

Finally we stress the following point: the potentiality of Fourier estimation method relies on the fact that it allows to reconstruct the volatility as a *stochastic function of time* in the univariate and multivariate case. This property makes possible to iterate the volatility functor and to compute the volatility of the volatility; see [33] for an early recognition on this point. In other words, we can handle the volatility function as an observable variable. For instance, this feature is essential when a stochastic derivation of volatility along the time evolution is performed as in contingent claim pricing-hedging (see [11, 35]), or when we study the geometry of the Heath–Jarrow–Morton interest rates dynamics given the observation of a single market trajectory [34].

The paper is organized as follows. Section 2 contains the main result: we derive the Fourier transform of the volatility function as the Bohr convolution product of the Fourier transform of log-returns. In Section 3, given a discrete, unevenly spaced and asynchronously sampled observations of the price processes, we prove the consistency uniformly in time of the random estimator of the co-volatility function by considering a fixed time span and increasing the frequency of the observations. In Section 4, the weak convergence in distribution of the co-volatility process to a mixture of Gaussian laws is proved. In Section 5 we analyze the implications of our results for the Fourier estimator of integrated volatility and co-volatility, even in the presence of microstructure noise effects. Section 6 concludes. Proofs are in Section 7.

**2. Instantaneous co-volatility computation.** Let $p(t) = (p^1(t), \ldots, p^n(t))$ be a continuous semi-martingale satisfying the following Itô stochastic differential equations

$$(2) \qquad dp^j(t) = \sum_{i=1}^{d} \sigma_i^j(t) \, dW^i + b^j(t) \, dt, \qquad j = 1, \ldots, n,$$

where $W = (W^1, \ldots, W^d)$ are independent Brownian motions on a filtered probability space satisfying the *usual conditions*, and $\sigma_*^*$ and $b^*$ are adapted random processes satisfying

$$(H) \qquad \begin{aligned} & E\left[\int_0^T (b^i(t))^2 \, dt\right] < \infty, \\ & E\left[\int_0^T (\sigma_i^j(t))^4 \, dt\right] < \infty, \qquad i = 1, \ldots, d, \; j = 1, \ldots, n. \end{aligned}$$

In financial applications the process $p(t)$ represents the log-price of assets and from the representation (2) we define the *volatility matrix*, which in our hypothesis *depends upon time*:

$$(3) \qquad \Sigma^{j,k}(t) = \sum_{i=1}^{d} \sigma_i^j(t) \sigma_i^k(t).$$

For notational simplicity we will refer to the case of two assets whose prices are $(p^1(t), p^2(t))$. We compute the so called *instantaneous volatilities* $\Sigma^{jj}(t)$ for $j = 1, 2$ and *co-volatilities* $\Sigma^{12}(t) = \Sigma^{21}(t)$.

By change of the origin of time and rescaling the unit of time we can always reduce ourselves to the case where the time window $[0, T]$ becomes $[0, 2\pi]$. We recall some definitions from harmonic analysis theory (see, e.g., [30]): given a function $\phi$ on the circle $S^1$, we consider its *Fourier transform*, which is defined on the group of integers $\mathbf{Z}$ by the formula

$$\mathcal{F}(\phi)(k) := \frac{1}{2\pi} \int_0^{2\pi} \phi(\vartheta) \exp(-\mathrm{i}k\vartheta) \, d\vartheta \qquad \text{for } k \in \mathbf{Z}.$$

We define

$$\mathcal{F}(d\phi)(k) := \frac{1}{2\pi} \int_{]0,2\pi[} \exp(-\mathrm{i}k\vartheta) \, d\phi(\vartheta),$$

and, using integration by parts, we remark that

$$\mathcal{F}(\phi)(k) = \frac{i}{k}\left[\frac{1}{2\pi}(\phi(2\pi) - \phi(0)) - \mathcal{F}(d\phi)(k)\right].$$

Given two functions $\Phi, \Psi$ on the integers $\mathbf{Z}$, we say that the *Bohr convolution product* exists if the following limit exists for all integers $k$

$$(\Phi *_B \Psi)(k) := \lim_{N \to \infty} \frac{1}{2N+1} \sum_{s=-N}^{N} \Phi(s) \Psi(k-s). \tag{4}$$

The subindex $B$ is the initial of *Bohr*, who developed similar ideas in the context of the theory of almost periodic functions.

The main result of this section is formula (5) on which the construction of the volatility estimator relies. Theorem 2.1 contains the identity relating the Fourier transform of the price process $p(t)$ to the Fourier transform of the volatility matrix $\Sigma(t)$.

THEOREM 2.1. *Consider a process $p$ satisfying the assumption (H). Then we have for $i, j = 1, 2$:*

$$\frac{1}{2\pi} \mathcal{F}(\Sigma^{ij}) = \mathcal{F}(dp^i) *_B \mathcal{F}(dp^j). \tag{5}$$

*The convergence of the convolution product (5) is attained in probability.*

The proof is based on the Itô energy identity for stochastic integrals. A preliminary step proves that the drift component $b^*(t)$ of the semi-martingale $p(t)$ gives no contribution to the formula (5). The latter derives from the following result.

LEMMA 2.2. *Consider a function $u$ defined on $\mathbf{Z}$ and define the Bohr $L^2$-seminorm*

$$\|u\|_{BL^2}^2 := \limsup_{N \to \infty} \frac{1}{2N+1} \sum_{k=-N}^{N} (u(k))^2.$$

*Then*

$$\|u *_B v\|_{L^\infty} \leq \|u\|_{BL^2} \|v\|_{BL^2}. \tag{6}$$

*In particular, it follows that*

$$\|v\|_{BL^2} = 0 \quad \Rightarrow \quad (u *_B v) = 0 \quad \text{for all } u. \tag{7}$$

By Theorem 2.1 we gather all the Fourier coefficients of the volatility matrix by means of the Fourier transform of the log-returns. Then the reconstruction of the cross-volatility functions $\Sigma^{ij}(t)$ from its Fourier coefficients, can be obtained as follows: define for $i, j = 1, 2$

$$\Phi_N^i(k) := \mathcal{F}(dp^i)(k) \quad \text{for } |k| \leq 2N \quad \text{and} \quad 0 \text{ otherwise},$$

and, for any $|k| \leq N$,

$$\Psi_N^{ij}(k) := \frac{1}{2N+1} \sum_{s \in \mathbf{Z}} \Phi_N^i(s) \Phi_N^j(k-s).$$

If $\Sigma^{ij}(t)$ is continuous, then the Fourier–Fejer summation gives almost everywhere

(8) $\quad \Sigma^{ij}(t) = \lim_{N \to \infty} \sum_{|k| < N} \left(1 - \frac{|k|}{N}\right) \Psi_N^{ij}(k) \exp(\mathrm{i}kt) \qquad$ for all $t \in (0, 2\pi)$.

REMARK 2.3. We note that in the case of the volatility functions $\Sigma^{ii}(t)$, the approximating trigonometric functions appearing in the sequence (8) are positive. In fact, by the positivity of the Fourier transform of the Fejer kernel it is enough to show the same result for the function

$$Q_N(t) := \sum_{k \in \mathbf{Z}} \Psi_N^i(k) \exp(\mathrm{i}kt).$$

Then we have

$$(2N+1) Q_N(t) = \sum_{k \in \mathbf{Z}} \exp(\mathrm{i}kt) \sum_{s \in \mathbf{Z}} \Phi_N^i(s) \Phi_N^i(k-s)$$

and this is equal to

$$\sum_{s, s' \in \mathbf{Z}} \exp(\mathrm{i}(s - s')t) \Phi_N^i(s) \overline{\Phi^i}_N(s') = \left| \sum_{s \in \mathbf{Z}} \exp(\mathrm{i}st) \Phi_N^i(s) \right|^2,$$

using the fact that $\overline{\Phi^i}_N(-s) = \Phi_N^i(s)$, where $\overline{\Phi^i}$ denotes the conjugate function of $\Phi^i$.

In the case of a single asset, the computation of the volatility function follows by Theorem 2.1. For completeness we state the result.

COROLLARY 2.4. *Suppose that the log-price process $p(t)$ follows the Itô process*

$$dp(t) = b(t)\, dt + \sigma(t)\, dW(t),$$

*where $W$ is a standard Brownian motion and $\sigma$, $b$ are adapted random processes satisfying assumption (H). Then the Fourier coefficients of the volatility function $\sigma^2(t)$ are obtained, for any integer $k$, by*

$$\frac{1}{2\pi} \mathcal{F}(\sigma^2)(k) = \lim_{N \to \infty} \frac{1}{2N+1} \sum_{|s| \leq N} \mathcal{F}(dp)(s) \mathcal{F}(dp)(k-s),$$

*where the limit is attained in probability.*

The instantaneous volatility function $\sigma^2(t)$ can be obtained analogously with (8), using the Fourier–Fejer summation.

**3. Consistency of the Fourier estimator of instantaneous co-volatility.** In this section we prove that the estimation methodology proposed does not induce any synchronization bias and the Fourier estimator of instantaneous co-volatility is consistent. As a particular case the consistency of the integrated co-volatility estimator follows.

Let $\mathcal{S}_{n_1} := \{t^1_{i,n_1}, i = 1, \ldots, n_1\}$ and $\mathcal{T}_{n_2} := \{t^2_{j,n_2}, j = 1, \ldots, n_2\}$ be the trades for the asset 1 and 2, respectively. For the ease of notation we often suppress the second index $n_1, n_2$. For simplicity suppose that both assets trade at $t_0 = 0$ and $t^1_{n_1} = t^2_{n_2} = 2\pi$. It is not restrictive to suppose that $n_1 = n_2 =: n$. Denote for $k = 1, 2$, $\rho^k(n) := \max_{0 \leq h \leq n-1} |t^k_{h+1} - t^k_h|$ and suppose that $\rho(n) := \rho^1(n) \vee \rho^2(n) \to 0$ as $n \to \infty$.

Consider the following interpolation formula

$$p^1_n(t) := \sum_{i=0}^{n-1} p^1(t^1_i) I_{[t^1_i, t^1_{i+1}[}(t).$$

Analogously for $p^2_n(t)$. We will indicate $I^1_i := [t^1_i, t^1_{i+1}[$ and $J^2_j := [t^2_j, t^2_{j+1}[$ and the returns by $\delta_{I^1_i}(p^1) := p^1(t^1_{i+1}) - p^1(t^1_i)$ and $\delta_{J^2_j}(p^2) := p^2(t^2_{j+1}) - p^2(t^2_j)$.

For any integer $k$, $|k| \leq 2N$, let

(9)
$$c_k(dp^1_n) := \frac{1}{2\pi} \sum_{i=0}^{n-1} \exp(-\mathrm{i}kt^1_i) \delta_{I^1_i}(p^1),$$

$$c_k(dp^2_n) := \frac{1}{2\pi} \sum_{j=0}^{n-1} \exp(-\mathrm{i}kt^2_j) \delta_{J^2_j}(p^2),$$

moreover define

(10) $$c_k(\Sigma^{12}) := \frac{1}{2\pi} \int_0^{2\pi} e^{-\mathrm{i}kt} \Sigma^{12}(t) \, dt.$$

Let for any $|k| \leq N$

(11) $$\alpha_k(N, p^1_n, p^2_n) := \frac{2\pi}{2N+1} \sum_{|s| \leq N} c_s(dp^1_n) c_{k-s}(dp^2_n).$$

Finally, define

(12) $$\widehat{\Sigma}^{12}_{n,N}(t) := \sum_{|k| \leq N} \left(1 - \frac{|k|}{N}\right) \alpha_k(N, p^1_n, p^2_n) e^{\mathrm{i}kt}.$$

The random function $\widehat{\Sigma}^{12}_{n,N}(t)$ will be called the *Fourier estimator* of the instantaneous co-volatility $\Sigma^{12}(t)$. In the following theorems we prove the convergence in probability of the random function (12) to the co-volatility function $\Sigma^{12}(t)$ uniformly in $t$, as $N$ and $n$ go to infinity.

Before stating the result, some comments are needed. We note that the effective reconstruction of $\Sigma^{12}(t)$ is realized as superposition of three limits. The first limiting procedure consists in the convergence as $N$ goes to infinity of the Fejer sum

$$\sum_{|k|<N}\left(1-\frac{|k|}{N}\right)\exp(ikt)c_k(\Sigma^{12}).$$

To this aim we remark that, given a continuous function $\phi$, denote by $\omega_\phi(\lambda)$ its modulus of continuity defined as

$$\omega_\phi(\lambda) := \sup_{|\theta-\theta'|\leq\lambda} |\phi(\theta)-\phi(\theta')|$$

and let $c_k(\phi) := \mathcal{F}(\phi)(k)$; then it holds

(13) $$\sup_t \left|\phi(t) - \sum_{|k|<N}\left(1-\frac{|k|}{N}\right)\exp(ikt)c_k(\phi)\right| \leq \omega_\phi\left(\frac{4}{N}\right)$$

(see [40] for the converse of this inequality). Therefore, if we assume that some a priori information on the modulus of continuity of $\Sigma(t)$ is available, from (13) we see that the determinacy of $\Sigma^{12}(t)$ in the $L^\infty$-norm can be obtained from the knowledge of a finite number of its Fourier coefficients. In the sequel suppose that $\Sigma^{ij}(t)$ satisfies

(H$_1$)    $\operatorname{ess\,sup}\|\Sigma^{ij}\|_{L^\infty} < \infty,$    where $\|\Sigma^{ij}\|_{L^\infty} := \sup_t |\Sigma^{ij}(t)|.$

The determination of $c_k(\Sigma^{12})$ depends upon the superposition of two limiting procedures, letting $n$ go to infinity. First, we approximate $c_k(dp^1)$ by $c_k(dp_n^1)$. The second limit procedure appears in the definition of the Bohr convolution. Next, two lemmas provide the required inequalities to obtain the consistency result.

LEMMA 3.1. *For any integers $k,s$, with $|k|\leq N$ and $|s|\leq N$, and any $n$, the following inequality holds:*

(14) $$E[|c_s(dp_n^1)c_{k-s}(dp_n^2) - c_s(dp^1)c_{k-s}(dp^2)|^2]$$
$$\leq CN^2\rho(n)^2 \operatorname{ess\,sup}\|\Sigma^{11}\|_{L^\infty} \operatorname{ess\,sup}\|\Sigma^{22}\|_{L^\infty},$$

*where $C$ is a suitable constant.*

LEMMA 3.2. *For any $k$, $|k|\leq 2N$, the following inequality holds:*

(15) $$E\left[\left|\frac{2\pi}{2N+1}\sum_{|s|\leq N} c_s(dp^1)c_{k-s}(dp^2) - c_k(\Sigma^{12})\right|^2\right]$$
$$\leq \frac{2}{2N+1}\operatorname{ess\,sup}\|\Sigma^{11}\|_{L^\infty}\operatorname{ess\,sup}\|\Sigma^{22}\|_{L^\infty}.$$

THEOREM 3.3. *Given $c_k(dp_n^1)$ and $c_k(dp_n^2)$ defined in (9) and $c_k(\Sigma^{12})$ defined in (10), let for any $|k| \leq N$*

$$\text{(16)} \quad \alpha_k(N, p_n^1, p_n^2) := \frac{2\pi}{2N+1} \sum_{|s| \leq N} c_s(dp_n^1) c_{k-s}(dp_n^2).$$

*Suppose that $\Sigma(t)$ is continuous and $N\rho(n) \to 0$ as $N, n \to \infty$. Then, for any $k$, the following convergence in probability holds:*

$$\text{(17)} \quad \lim_{n,N \to \infty} \alpha_k(N, p_n^1, p_n^2) = c_k(\Sigma^{12}).$$

Finally, we conclude by [40]:

THEOREM 3.4. *Let $\widehat{\Sigma}^{12}_{n,N}(t)$ defined in (12). If $\Sigma(t)$ is continuous and $N\rho(n) \to 0$, the following convergence holds in probability*

$$\lim_{n,N \to \infty} \sup_{0 \leq t \leq 2\pi} |\widehat{\Sigma}^{12}_{n,N}(t) - \Sigma^{12}(t)| = 0.$$

In the univariate model specified in Corollary 2.4 the Fourier estimator of instantaneous volatility can be obtained as follows.

Consider a discrete unevenly spaced sampling of the process $p$. Fix a sequence $\mathcal{S}_n$ of finite subsets of $[0, 2\pi]$, let $\mathcal{S}_n := \{0 = t_{0,n} \leq t_{1,n} \leq \cdots \leq t_{k_n,n} = 2\pi\}$ for any $n \geq 1$ such that $\rho(n) := \max_{0 \leq h \leq k_n - 1} |t_{h+1,n} - t_{h,n}| \to 0$ as $n \to \infty$. Denote, for any $j$, $\delta_j(p) := p(t_{j+1,n}) - p(t_{j,n})$.

Define, for $|k| \leq N$,

$$\text{(18)} \quad \alpha_k(N, p_n) := \frac{2\pi}{2N+1} \sum_{|s| \leq N} c_s(dp_n) c_{k-s}(dp_n),$$

where for any integer $k$, $|k| \leq 2N$, $c_k(dp_n)$ are the Fourier coefficients of the log-return process, that is,

$$c_k(dp_n) = \frac{1}{2\pi} \sum_{j=0}^{k_n - 1} \exp(-\mathrm{i}kt_{j,n}) \delta_j(p).$$

COROLLARY 3.5. *Consider for any $n, N$ the random function of time*

$$\text{(19)} \quad \widehat{\sigma}^2_{n,N}(t) := \sum_{|k| < N} \left(1 - \frac{|k|}{N}\right) \alpha_k(N, p_n) \exp(\mathrm{i}kt),$$

*where $\alpha_k(N, p_n)$ is defined in (18). Suppose $\sigma^2(t)$ is continuous.*

(i) *For any integer $k$ the following convergence in probability holds*

$$\lim_{n,N \to \infty} \alpha_k(N, p_n) = c_k(\sigma^2).$$

(ii) *It holds in probability*

$$\lim_{n,N\to\infty} \sup_{t\in[0,2\pi]} |\widehat{\sigma}^2_{n,N}(t) - \sigma^2(t)| = 0.$$

REMARK 3.6. In the construction of the Fourier estimator we can consider observations $\{t_h^k\}$ which are random but independent of the observed processes $p^k$, $k = 1, 2$. In fact by splitting the probability space, this kind of sampling reduces to a deterministic sampling. Therefore all our study has been made in this context.

**4. Limiting distribution.** In this section we prove the weak convergence of the instantaneous co-volatility estimator defined in (12) to the process $\Sigma^{12}(t)$. As a corollary the asymptotic normality for the error of the integrated co-volatility estimator follows. The asymptotic results are obtained as the mesh of the partition $\rho(n) \to 0$ and the interval $[0, 2\pi]$ remains fixed. Let $N$ be the number of Fourier coefficients considered in the reconstruction of $\Sigma^{ij}(t)$ and suppose that $N \to \infty$.

We first consider the univariate model specified in Corollary 2.4 and the estimator of the volatility function $\widehat{\sigma}^2_{n,N}(t)$ given by (19). As we suppose to have an irregularly spaced grid, some conditions are needed. The following assumptions from [37] are sufficient:

(A)  (i) $\rho(n) \to 0$ and $k_n \rho(n) = 0(1)$,
   (ii) $H_n(t) := \dfrac{\sum_{t_{j+1,n} \leq t}(t_{j+1,n} - t_{j,n})^2}{2\pi/k_n} \to H(t)$ as $n \to \infty$,
   (iii) $H(t)$ is continuously differentiable.

In the sequel we will denote by $\overline{\psi}$ the conjugate of the complex function $\psi$.

THEOREM 4.1. *Let $\widehat{\sigma}^2_{n,N}(t)$ defined in (19). Assume that as $n \to \infty$ then $\rho(n) \to 0$, $N \to \infty$, and assumption* (A) *is satisfied. Then, for any function $h \in \text{Lip}(\alpha)$, $\alpha > \frac{1}{2}$, with compact support in $(0, 2\pi)$,*

$$(\rho(n))^{-1/2} \int_0^{2\pi} h(t)(\overline{\widehat{\sigma}^2_{n,N}(t)} - \sigma^2(t))\, dt$$

*converges in law to a mixture of Gaussian distribution with random variance $2\int_0^{2\pi} H'(t) h^2(t) \sigma^4(t)\, dt$, provided that $\rho(n) N^{2\alpha} \to \infty$.*

The result in Theorem 4.1 derives from a central limit result of [27], generalized for unevenly spaced data in [37] (Proposition 1), after the following representation theorem holds.

THEOREM 4.2. *Under the hypothesis of Theorem 4.1 the following result holds*

$$(\rho(n))^{-1/2} \int_0^{2\pi} h(t)(\overline{\widehat{\sigma}_{n,N}^2(t)} - \sigma^2(t)) \, dt$$

(20)
$$= (\rho(n))^{-1/2} 2 \sum_{j=0}^{k_n-1} h_N(t_{j,n}) \int_{t_{j,n}}^{t_{j+1,n}} \left( \int_{t_{j,n}}^t \sigma(s) \, dW(s) \right) \sigma(t) \, dW(t)$$

$$+ o_p(1),$$

*where the function $h_N$ is defined as $h_N(t) := \sum_{|k|<N}(1 - \frac{|k|}{N}) c_k(h) e^{ikt}$ and $o_p(1)$ is a term which goes to zero in probability.*

REMARK 4.3. The above study of the asymptotic properties of the estimator also provides a way to obtain the order of magnitude of the number of Fourier coefficients, $N$, to be included in the estimation given the number of observations, $n$. The asymptotic properties require that $\rho(n)N \to 0$ and $\rho(n)N^{2\alpha} \to \infty$, if $\rho(n)$ is the mesh of the partition. Similar conditions appear in order to derive the limiting properties of Bartlett type kernel estimator (see [43]) for the relative growth of the number of auto-covariances with respect to the number of data.

The asymptotic result in the bivariate setting needs to take into account for the asynchronicity of the returns. Then assumption (A)(ii) must be modified according to

(A) (ii') $H_n(t) := \frac{n}{2\pi} \sum_{t_{i+1}^1 \wedge t_{j+1}^2 \leq t} (t_{i+1}^1 \wedge t_{j+1}^2 - t_i^1 \vee t_j^2)^2 I_{\{t_i^1 \vee t_j^2 < t_{i+1}^1 \wedge t_{j+1}^2\}} \to H(t)$ as $n \to \infty$.

THEOREM 4.4. *Let $\widehat{\Sigma}_{n,N}^{12}(t)$ defined in (12). Suppose $\rho(n) := \rho^1(n) \vee \rho^2(n) \to 0$, $N \to \infty$, and assumption (A) (i), (ii'), (iii) hold, then, for any function $h \in \text{Lip}(\alpha)$, $\alpha > \frac{2}{3}$, with compact support in $(0, 2\pi)$,*

(21) $$(\rho(n))^{-1/2} \int_0^{2\pi} h(t)(\overline{\widehat{\Sigma}_{n,N}^{12}(t)} - \Sigma^{12}(t)) \, dt$$

*converges in law to a mixture of Gaussian distribution with variance $\int_0^{2\pi} H'(t) \times h^2(t)(\Sigma^{11}(t)\Sigma^{22}(t) + (\Sigma^{12}(t))^2) \, dt$, provided that $\rho(n)N^{2\alpha} \to \infty$, if $\alpha > \frac{2}{3}$, and $\rho(n)N^{4/3} \to 0$.*

In order to apply central limit result of [27] (see also Section 4 in [9]), the following representation theorem is required.

A FOURIER TRANSFORM METHOD FOR NONPARAMETRIC ESTIMATION 13THEOREM 4.5. *Under the hypothesis of Theorem 4.4, the following result holds*

$$(\rho(n))^{-1/2} \int_0^{2\pi} h(t)(\overline{\widehat{\Sigma}_{n,N}^{12}(t)} - \Sigma^{12}(t))\, dt$$

$$= (\rho(n))^{-1/2} \sum_{i,j} h_N(t_j^2) \left( \int_{t_i^1 \vee t_j^2}^{t_{i+1}^1 \wedge t_{j+1}^2} \int_{t_i^1 \vee t_j^2}^{t} dp^1(s)\, dp^2(t) \right.$$

$$\left. + \int_{t_i^1 \vee t_j^2}^{t_{i+1}^1 \wedge t_{j+1}^2} \int_{t_i^1 \vee t_j^2}^{t} dp^2(s)\, dp^1(t) \right) + o_p(1),$$

where the function $h_N$ is defined as $h_N(t) := \sum_{|k|<N}(1 - \frac{|k|}{N})c_k(h)e^{ikt}$ and $o_p(1)$ is a term which goes to zero in probability.

**5. Fourier estimator of integrated volatilities in the presence of microstructure noise.** The financial econometric's literature mainly focus on the *integrated* volatility and co-volatility, that is, $\int_0^T \Sigma^{ij}(t)\, dt$, for $i, j = 1, 2$, where $[0, T]$ is a fixed time horizon, for example, a day. In fact the computation of spot volatility by using quadratic variation formula is quite instable, because it involves a numerical derivative.

In the context of Fourier estimation methodology the integrated volatility and co-volatility are computed by considering the $k = 0$ coefficient, respectively, in the formulae (11) and (18). The theory presented in this paper is stronger because we recover pathwise the entire volatility curve from the observation of an asset price trajectory. Nonetheless, due to the importance of the integrated (co-) volatility for financial econometrics purposes, in this section we explicitly state the results concerning integrated volatilities; in particular, the properties of Fourier estimator in comparison with different integrated volatility estimators proposed in the literature can be analyzed in this context.

Consider first the univariate case. Suppose that the price process $p(t)$ is a continuous semi-martingale satisfying the hypothesis stated in Corollary 2.4. Observe that by definition

$$2\pi \mathcal{F}(\sigma^2)(0) = \int_0^{2\pi} \sigma^2(t)\, dt,$$

therefore the integrated volatility over a given period (always reducible to $[0, 2\pi]$) is computed by Corollary 2.4

(22) $$\int_0^{2\pi} \sigma^2(t)\, dt = P - \lim_{N \to \infty} \frac{(2\pi)^2}{2N+1} \sum_{|s| \leq N} \mathcal{F}(dp)(s)\mathcal{F}(dp)(-s).$$

More precisely, using the same notation as in Section 3, from (22) the *Fourier estimator of integrated volatility* over $[0, 2\pi]$ is defined as follows

$$
(23) \qquad \hat{\sigma}^2_{n,N} := \frac{(2\pi)^2}{2N+1} \sum_{|s| \leq N} c_s(dp_n) c_{-s}(dp_n).
$$

The consistency in probability uniformly in time of the Fourier estimator of the integrated volatility follows immediately by Corollary 3.5. The efficiency of the Fourier procedure to obtain estimator of integrated volatility is studied in [13, 21]. We remark that the Fourier estimator (23) can be expressed as

$$
(24) \qquad \hat{\sigma}^2_{n,N} = \sum_{j=1}^{k_n} \sum_{j'=1}^{k_n} D_N(t_j - t'_j) \delta_j(p) \delta_{j'}(p),
$$

where $D_N(t)$ denotes the rescaled Dirichlet kernel defined by

$$
(25) \qquad D_N(t) := \frac{1}{2N+1} \sum_{|s| \leq N} e^{\mathrm{i}st} = \frac{1}{2N+1} \frac{\sin[(N+1/2)t]}{\sin(t/2)}.
$$

From (24) different features of the Fourier estimation method are highlighted. First, the Fourier estimator uses all available data by integration, thus incorporating not only the squared log-returns but also the products of disjoint log-returns along the time window. Second, the convolution product weights the auto-covariances at any given frequency. In fact it is easily seen from (24) that

$$
\hat{\sigma}^2_{n,N} = RV + \sum_{j \neq j'} D_N(t_j - t'_j) \delta_j(p) \delta_{j'}(p),
$$

where $RV$ denotes the realized volatility estimator, namely $\sum_{j=1}^{k_n} (\delta_j(p))^2$. The contribution of various order auto-covariances has early been considered by [44] and recently used to correct the bias of the realized variance type estimators in the presence of microstructure noise, in particular the subsampled estimator by [43] and the realized (subsampled) kernels by [10]. Different weighting functions are also considered in [6, 29] who express many volatility estimators proposed in the literature in an unified way as a quadratic form. In the analysis of [23, 29] a peculiar feature of Fourier estimator is not efficiently addressed: the Dirichlet kernel appearing in the Fourier volatility estimator depends on the number of frequencies $N$, besides the number of auto-covariances and the time lag between two observations. The importance of the choice of the frequencies has been empirically noticed by [38]; the authors make an empirical comparison, through Monte Carlo simulations and using high frequency market data, between some estimators of integrated volatility. Their analysis shows that the Fourier estimator is superior to the realized volatility and the wavelet estimators and, even compared to some bias correcting methods for microstructure noise, namely the kernel-based estimator by [44] and the related unbiased estimator proposed by [22], the Fourier estimator provides smaller root mean squared error while having only slightly higher bias. Recently [36] prove that, while the MSE of the realized volatility estimator under microstructure noise diverges as the number $n$ of observations increases, the MSE of the Fourier estimator is substantially unaffected by the presence of microstructure noise by choosing in an appropriate way the number of Fourier coefficients to be included in the estimation, as indicated explicitly by the MSE computation. The optimal MSE-based Fourier estimator turns out to be very efficient, even in comparison with methods specifically designed to handle market microstructure contaminations; more specifically, the Fourier estimator is competitive in terms of MSE for high sampling frequencies up to 30 seconds. We remark that in the presence of microstructure noise it is useful to consider the following more stable version of Fourier estimator

$$(26) \qquad \widetilde{\sigma}^2_{n,N} := \frac{(2\pi)^2}{N+1} \sum_{s=-N}^{N} \left(1 - \frac{|s|}{N}\right) c_s(dp_n) c_{-s}(dp_n).$$

As for the multivariate case, the integrated co-volatility is recovered as a byproduct of formula (5) through the following identity

$$2\pi \mathcal{F}(\Sigma^{12})(0) = \int_0^{2\pi} \Sigma^{12}(t)\, dt.$$

Therefore the *Fourier estimator of integrated co-volatility* is defined as

$$(27) \qquad \widehat{\Sigma}^{12}_{n,N} := \frac{1}{2N+1} \sum_{|s|\leq N} \sum_{i=0}^{n-1} \sum_{j=0}^{n-1} e^{(is(t_i^1 - t_j^2))} \delta_{I_i^1}(p^1) \delta_{J_j^2}(p^2).$$

The consistency of the Fourier estimator of integrated co-volatility with asynchronous observations for different assets follows from Theorem 3.3. Moreover by (14) and (15) we get an estimate of the mean squared error (MSE) of the Fourier integrated co-volatility estimator

$$MSE_{\widehat{\Sigma}^{12}_{n,N}} \leq \left(C_1 \frac{1}{N} + C_2 N^2 \rho(n)^2\right) \operatorname{ess\,sup} \|\Sigma^{11}\|_{L^\infty} \operatorname{ess\,sup} \|\Sigma^{22}\|_{L^\infty},$$

$C_1$ and $C_2$ being constants. By choosing $N$ such that $\frac{1}{N} = N^2 \rho(n)^2$ (up to a constant) we obtain that $N = \rho(n)^{-2/3}$. Therefore the mean squared error of Fourier estimator of co-volatility is of order $O(\rho(n)^{2/3})$.

All the realized covariance type estimators rely on "synchronization" methods; therefore, they suffer from a downward bias when the sampling interval is reduced (known as *Epps effect* by [16]). Recently [20] have proposed an estimator of integrated co-volatility which is consistent under asynchronous observations. This estimator is defined by adding the cross-products of all fully and partially overlapping returns. [39] prove that the estimator of [20] is biased in the presence of microstructure noise and its variance diverges as the number of observations increases. We note that Fourier estimator of integrated co-volatility defined in (27) has the advantage that from one side it is consistent in the absence of microstructure contaminations, but with the optimal choice of $N$ (as explained in the univariate case), it is more robust in the presence of microstructure noise effects.

In summary, the Fourier estimator exploits more information in the data without being affected by a severe bias coming from microstructure contamination, simply by choosing the lower frequencies.

**6. Conclusions.** In this paper we have proposed a nonparametric estimator of the instantaneous volatility and co-volatility function in the context of continuous semi-martingale models, given the knowledge of the Fourier transform of the log-returns. The method applies both to univariate and multivariate setting and allows to consider unevenly spaced and asynchronous observations. Moreover we have studied the statistical properties of the Fourier estimator by providing consistency and asymptotic normality. We emphasize two main features of the proposed estimation methodology: the first one is that the definition of Fourier estimator is designed specifically to handle the computation of co-volatilities without any manipulation of data, due to the fact that the Fourier method requires an integration of all data at disposal for any asset; the second one is that, with the availability of high frequency data, we are able to reconstruct the co-volatility as a stochastic function in an effective way, which permits to handle the volatility function as an observable variable in the econometric and financial applications.

**7. Proofs.** In the sequel we will denote by $\overline{\psi}$ the conjugate of the complex function $\psi$.

PROOF OF LEMMA 2.2. For $a \in \mathbf{Z}$ denote $u^a(k) := u(k+a)$ and $\check{u}(k) := \overline{u(-k)}$. Observe that

$$(u *_B v)(a) = (u|\check{v}^{-a})_{BL^2},$$

where $(u|v)_{BL^2}$ denotes the scalar product

$$(u|v)_{BL^2} := \lim_{N \to \infty} \frac{1}{2N+1} \sum_{k=-N}^{N} u(k)\overline{v(k)}.$$

By Cauchy–Schwarz inequality it holds:

$$\frac{1}{2N+1} \sum_{k=-N}^{N} u(k)\overline{v(k)} \leq \left(\frac{1}{2N+1} \sum_{k=-N}^{N} (u(k))^2\right)^{1/2} \left(\frac{1}{2N+1} \sum_{k=-N}^{N} (v(k))^2\right)^{1/2}.$$

This implies that for all $a$
$$(u *_B v)(a) \leq \|u\|_{BL^2} \|\check{v}^{-a}\|_{BL^2}.$$

Because $\|\check{v}^{-a}\|_{BL^2} = \|v\|_{BL^2}$, the result is proved. □

PROOF OF THEOREM 2.1. Without loss of generality we suppose that $b^* = 0$. In fact, consider the process $p_m = (p_m^1, p_m^2)$ given by $dp_m^j(t) = \sum_{i=1}^{2} \sigma_i^j(t) \, dW^i(t)$: $p_m$ has the same volatility function than the process $p$. We have that for $i = 1, 2$:
$$\mathcal{F}(dp^j)(k) = \mathcal{F}(dp_m^j)(k) + \mathcal{F}(b^j)(k).$$

For notational simplicity we omit now the upper index $j = 1, 2$. With the notation $\Phi_m(k) := \mathcal{F}(dp_m)(k)$, $\Phi_b(k) := \mathcal{F}(b)(k)$, we now prove that
$$(\Phi_m + \Phi_b) *_B (\Phi_m + \Phi_b) = \Phi_m *_B \Phi_m,$$

because the terms $\Phi_b *_B \Phi_b$, $\Phi_b *_B \Phi_m$ and $\Phi_m *_B \Phi_b$ give zero contribution. First, we consider $\Phi_b *_B \Phi_b$. We have
$$(\Phi_b *_B \Phi_b)(k) = \lim_{N \to \infty} \frac{1}{2N+1} \sum_{s=-N}^{N} \Phi_b(s) \Phi_b(k-s).$$

It follows by the Lemma 2.2 that:
$$\|\Phi_b *_B \Phi_b\|_{L^\infty} \leq \|\Phi_b\|_{BL^2} \|\Phi_b\|_{BL^2}. \tag{28}$$

Moreover by Plancherel theorem it holds:
$$\sum_{k=-\infty}^{+\infty} |\mathcal{F}(b)(k)|^2 = \int_0^{2\pi} \frac{1}{2\pi} |b(t)|^2 \, dt. \tag{29}$$

Therefore by (29) it follows:
$$\|\Phi_b\|_{BL^2}^2 = \lim_{N \to \infty} \frac{1}{2N+1} \sum_{k=-N}^{N} |\mathcal{F}(b)(k)|^2 = 0. \tag{30}$$

The fact that the terms $\Phi_b *_B \Phi_m$ and $\Phi_m *_B \Phi_b$ give zero contribution follows by the Lemma 2.2 and (30).

Under the hypothesis that $b^* = 0$, then $p$ is a martingale. Introduce the complex martingales for any integer $k$ and $i = 1, 2$
$$\Gamma_k^i(t) := \frac{1}{2\pi} \int_0^t \exp(-iks) \, dp^i(s).$$

Then by definition $\Gamma_k^i(2\pi) = \mathcal{F}(dp^i)(k)$. Compute by Itô formula the stochastic differential
$$d(\Gamma_k^i \Gamma_r^j)(t) = \Sigma^{ij}(t) \exp(-i(k+r)t) \, dt + \Gamma_k^i(t) \, d\Gamma_r^j(t) + \Gamma_r^j(t) \, d\Gamma_k^i(t),$$

we obtain

(31) $$\Gamma^i_k(2\pi)\Gamma^j_r(2\pi) = \frac{1}{2\pi}\mathcal{F}(\Sigma^{ij})(k+r) + R^{ij}(k,r),$$

where

$$R^{ij}(k,r) := \int_0^{2\pi} \Gamma^i_k(t)\,d\Gamma^j_r(t) + \Gamma^i_r(t)\,d\Gamma^j_k(t).$$

Fix an integer $N \geq 1$ and define for any integer $q$ with $|q| \leq N$,

(32) $$\gamma^{ij}_q(N) = \frac{1}{2N+1}\sum_{s=-N}^{N}\Gamma^i_{q+s}(2\pi)\Gamma^j_{-s}(2\pi).$$

By definition (4), it results that, for any integer $q$

$$\lim_{N\to\infty}\gamma^{ij}_q(N) = (\mathcal{F}(dp^i) *_B \mathcal{F}(dp^j))(q).$$

Moreover by (31), we have:

(33) $$\gamma^{ij}_q(N) = \frac{1}{2\pi}\mathcal{F}(\Sigma^{ij})(q) + R^{ij}_N,$$

where

$$R^{ij}_N = \frac{1}{2N+1}\int_0^{2\pi}\sum_{l=-N}^{N}\Gamma^i_{q-l}(t)\,d\Gamma^j_l(t) + \Gamma^j_l(t)\,d\Gamma^i_{q-l}(t).$$

We prove that $R^{ij}_N$ converges to 0 in probability as $N \to \infty$. By symmetry we are reduced to study

$$A_N := \frac{1}{2N+1}\sum_{l=-N}^{N}\int_0^{2\pi}d\Gamma^j_l(t_2)\int_0^{t_2}d\Gamma^i_{q-l}(t_1).$$

Introducing the rescaled Dirichlet kernel

$$D_N(t) := \frac{1}{2N+1}\sum_{s=-N}^{N}\exp(ist) = \frac{1}{2N+1}\frac{\sin(N+1/2)t}{\sin(t/2)},$$

we have that

(34) $$A_N = \frac{1}{4\pi^2}\int_0^{2\pi}dp^j(t_2)\left[\int_0^{t_2}\exp(iqt_1)D_N(t_1-t_2)\,dp^i(t_1)\right].$$

Define now

$$\alpha(t_2) := \int_0^{t_2}\cos(qt_1)D_N(t_1-t_2)\sum_{k=1}^{2}\sigma^i_k(t_1)\,dW^k(t_1),$$

$$\beta(t_2) := \int_0^{t_2}\sin(qt_1)D_N(t_1-t_2)\sum_{k=1}^{2}\sigma^i_k(t_1)\,dW^k(t_1).$$

By Itô energy identity for real stochastic integrals we have (e.g., [31])

$$16\pi^4 E[|A_N|^2] = \sum_{k=1}^{2} E\left[\int_0^{2\pi} (\alpha^2(t_2) + \beta^2(t_2))(\sigma_k^j(t_2))^2 \, dt_2\right].$$

By Cauchy–Schwarz inequality we have

$(16\pi^4 E[|A_N|^2])^2$

$$\leq 4 \sum_{k=1}^{2} E\left[\int_0^{2\pi} (\sigma_k^j(t_2))^4 \, dt_2\right]$$

$$\times \left\{ E\left[\int_0^{2\pi} \left(\int_0^{t_2} D_N(t_1 - t_2)\cos(qt_1) \sum_{k=1}^{2} \sigma_k^i(t_1) \, dW^k(t_1)\right)^4 dt_2\right]\right.$$

$$\left. + E\left[\int_0^{2\pi} \left(\int_0^{t_2} D_N(t_1 - t_2)\sin(qt_1) \sum_{k=1}^{2} \sigma_k^i(t_1) \, dW^k(t_1)\right)^4 dt_2\right]\right\}.$$

Because $E[\int_0^{2\pi} (\sigma_k^j(t))^4 \, dt] < \infty$ by hypothesis, we evaluate the other term using Burkholder–Gundy inequality

$$E\left[\int_0^{2\pi} \left(\int_0^{t_2} D_N(t_1 - t_2)\cos(qt_1) \sum_{k=1}^{2} \sigma_k^i(t_1) \, dW^k(t_1)\right)^4 dt_2\right]$$

$$\leq 4 E\left[\int_0^{2\pi} \int_0^{t_2} D_N^4(t_1 - t_2) \sum_{k=1}^{2} (\sigma_k^i(t_1))^4 \, dt_1 \, dt_2\right];$$

then making a change of variable $t_1 = u$, $t_2 - t_1 = v$, we get

$$= 8 E\left[\int_0^{2\pi} \sum_{k=1}^{2} (\sigma_k^i(u))^4 \, du\right] \int_0^{2\pi} D_N^4(v) \, dv.$$

We remark that $|D_N(v)| \leq 1$, therefore

$$\int_0^{2\pi} D_N^4(v) \, dv \leq \int_0^{2\pi} D_N^2(v) \, dv = \frac{2\pi}{2N+1}.$$

The last equality is a consequence of Plancherel equality. Therefore letting $N \to \infty$ from (33) we get the convergence in probability. □

PROOF OF LEMMA 3.1. For any $i = 1, 2$ it holds

$$c_k(dp_n^i) - c_k(dp^i) = \frac{1}{2\pi} \int_0^{2\pi} \beta(t) \, dp^i(t),$$

where

$$\beta(t) := \sum_{l=0}^{n-1} e^{ikt_l^i}(1 - e^{ik(t - t_l^i)}) I_{[t_l^i, t_{l+1}^i[}(t).$$

Now by Burkholder–Gundy inequality we have

(35) $$E[|c_k(dp_n^i) - c_k(dp^i)|^4] \leq Ck^4 \rho(n)^4 \operatorname{ess\,sup} \|\Sigma^{ii}\|_{L^\infty}^2,$$

where $C$ is a suitable constant. Therefore combining (35) and the inequality

$$E[|c_s(dp_n^1)c_{k-s}(dp_n^2) - c_s(dp^1)c_{k-s}(dp^2)|^2]$$
$$\leq 2(E[|c_s(dp^1)|^4]^{1/2} E[|c_{k-s}(dp_n^2) - c_{k-s}(dp^2)|^4]^{1/2}$$
$$+ E[|c_{k-s}(dp^2)|^4]^{1/2} E[|c_s(dp_n^1) - c_s(dp^1)|^4]^{1/2}),$$

we complete the proof. □

PROOF OF LEMMA 3.2. By Itô formula

$$E\left[\left|\frac{2\pi}{2N+1} \sum_{|s|<N} c_s(dp^1)c_{k-s}(dp^2) - c_k(\Sigma^{12})\right|^2\right]$$
$$\leq 2\left(E\left[\left|\frac{1}{2\pi}\int_0^{2\pi} dp^1(t_2) \int_0^{t_2} e^{ikt_1} D_N(t_2 - t_1)\,dp^2(t_1)\right|^2\right]\right.$$
$$\left. + E\left[\left|\frac{1}{2\pi}\int_0^{2\pi} dp^2(t_2) \int_0^{t_2} e^{ikt_1} D_N(t_2 - t_1)\,dp^1(t_1)\right|^2\right]\right),$$

where $D_N(t)$ is the rescaled Dirichlet kernel. Let

$$A_N(t) = \int_0^t e^{ikv} D_N(t - v)\,dp^2(v),$$

then by Itô energy identity and Hölder inequality

$$E\left[\left|\int_0^{2\pi} dp^1(t_2) A_N(t_2)\right|^2\right] = E\left[\int_0^{2\pi} A_N^2(t_2) \Sigma^{11}(t_2)\,dt_2\right]$$
$$\leq \frac{2\pi}{2N+1} \operatorname{ess\,sup}\|\Sigma^{11}\|_{L^\infty} \operatorname{ess\,sup}\|\Sigma^{22}\|_{L^\infty}$$

as $\int_0^{2\pi} D_N^2(v)\,dv = \frac{2\pi}{2N+1}$. □

PROOF OF THEOREM 3.3. Using the inequalities obtained in the Lemmas 3.1 and 3.2 it follows

$$E\left[\left|\frac{2\pi}{2N+1} \sum_{|s|\leq N} c_s(dp_n^1)c_{k-s}(dp_n^2) - c_k(\Sigma^{12})\right|^2\right]$$
$$\leq \left(C_1 \frac{1}{2N+1} + C_2 N^2 \rho(n)^2\right) \operatorname{ess\,sup}\|\Sigma^{11}\|_{L^\infty} \operatorname{ess\,sup}\|\Sigma^{22}\|_{L^\infty},$$

where $C_1$ and $C_2$ are constants. □

PROOF OF THEOREM 3.4. In the previous theorem we have proved for any fixed $k$ the convergence in probability of each coefficient $\alpha_k(N, p_n^1, p_n^2)$ to the Fourier coefficient $c_k(\Sigma^{12})$, as $n$ and $N$ go to $\infty$. Therefore the thesis follows by applying (13). $\square$

PROOF OF THEOREM 4.2. By using the Fourier series decomposition we have

$$\frac{1}{2\pi}\int_0^{2\pi} h(t)(\overline{\widehat{\sigma}_n^2(t)} - \sigma^2(t))\,dt$$
$$= \sum_{|k|<N}\left(1 - \frac{|k|}{N}\right)c_k(h)\overline{\alpha_k}(N,p_n) - \sum_{|k|<N} c_k(h)\overline{c_k}(\sigma^2) - \sum_{|k|\geq N} c_k(h)\overline{c_k}(\sigma^2),$$

where $c_k(f)$ denotes the $k$th Fourier coefficient of the function $f$, $\alpha_k(\cdot)$ is defined in (18). By Cauchy–Schwarz inequality and applying (13) we have

$$(36)\quad (\rho(n))^{-1/2}\sum_{|k|\geq N}\left(1 - \frac{|k|}{N}\right)c_k(h)\overline{c_k}(\sigma^2) \leq C(\rho(n))^{-1/2}N^{-\alpha}\|\sigma^2\|_{L^2},$$

where $C$ is a suitable constant. As $\rho(n)N^{2\alpha} \to \infty$, then (36) goes to 0 as $n$, $N \to \infty$.

Now we consider

$$2\pi\sum_{|k|<N}\left(1 - \frac{|k|}{N}\right)c_k(h)\overline{\alpha_k}(N,p_n)$$
$$= \sum_{j=0}^{k_n-1} h_N(t_{j,n})(\delta_j(p))^2 + 2\sum_{|k|<N}\left(1 - \frac{|k|}{N}\right)c_k(h)A_N(k),$$

where

$$(37)\qquad h_N(t) := \sum_{|k|<N}\left(1 - \frac{|k|}{N}\right)c_k(h)e^{ikt}$$

and

$$A_N(k) := \sum_{j<j'} e^{ikt_{j,n}}\frac{1}{2N+1}\sum_{|h|<N} e^{ih(t_{j',n}-t_{j,n})}\delta_j(p)\delta_{j'}(p).$$

It holds

$$E\left[\left|\sum_{|k|<N}\left(1 - \frac{|k|}{N}\right)c_k(h)A_N(k)\right|^2\right]$$
$$(38)\qquad = E\left[\left|\sum_{j<j'} h_N(t_{j,n})D_N(t_{j',n}-t_{j,n})\delta_j(p)\delta_{j'}(p)\right|^2\right]$$

$$\leq C_1 \operatorname{ess\,sup} \|\sigma^2\|_{L^\infty}^2 \frac{2\pi}{2N+1},$$

where $C_1$ is a constant. It follows that $(\rho(n))^{-1/2} E[|\sum_{|k|<N}(1-\frac{|k|}{N})c_k(h) \times A_N(k)|^2]$ goes to 0 in probability as $n, N \to \infty$. Therefore it holds

$$(\rho(n))^{-1/2} 2\pi \sum_{|k|<N} \left(1 - \frac{|k|}{N}\right) c_k(h) \overline{\alpha_k}(N, p_n)$$

(39)

$$= (\rho(n))^{-1/2} \sum_{j=0}^{k_n-1} h_N(t_{j,n})(\delta_j(p))^2 + o_p(1).$$

Finally it remains to consider

$$\sum_{j=0}^{k_n-1} h_N(t_{j,n})(\delta_j(p))^2 - \sum_{|k|<N} c_k(h) \overline{c_k}(\sigma^2).$$

It is enough to study

$$\sum_{j=0}^{k_n-1} h_N(t_{j,n})(\delta_j(p))^2 - \int_0^{2\pi} h_N(t) \sigma^2(t) \, dt$$

by Itô calculus this is equal to

$$2 \sum_{j=0}^{k_n-1} h_N(t_{j,n}) \left( \int_{t_{j,n}}^{t_{j+1,n}} (p(t) - p(t_{j,n})) \sigma(t) \, dW_t + \int_{t_{j,n}}^{t_{j+1,n}} \sigma^2(t) \, dt \right)$$
$$- \int_0^{2\pi} h_N(t) \sigma^2(t) \, dt.$$

Denote

$$F_N := \sum_{j=0}^{k_n-1} \int_{t_{j,n}}^{t_{j+1,n}} (h_N(t_{j,n}) - h_N(t)) \sigma^2(t) \, dt,$$

it holds

$$E[|F_N|] \leq \left\{ 2 \sup_t |h_N(t) - h(t)| + \omega_h(\rho(n)) \right\} E\left[ \int_0^{2\pi} \sigma^2(t) \, dt \right]$$
$$= \left\{ O\left(\frac{2}{N^\alpha}\right) + O(\rho(n)^\alpha) \right\} E\left[ \int_0^{2\pi} \sigma^2(t) \, dt \right].$$

Therefore by the condition $\rho(n) N^{2\alpha} \to \infty$ as $n \to \infty$ we have

$$(\rho(n))^{-1/2} E[|F_N|] \to 0 \quad \text{as } n \to \infty. \qquad \square$$

PROOF OF THEOREM 4.1. By the representation result in Theorem 4.2, we study the convergence in distribution of

$$(\rho(n))^{-1/2} 2 \sum_{j=0}^{k_n-1} h_N(t_{j,n}) \int_{t_{j,n}}^{t_{j+1,n}} \left( \int_{t_{j,n}}^{t} \sigma(s) \, dW_s \right) \sigma(t) \, dW_t.$$

We write

$$2 \sum_{j=0}^{k_n-1} h_N(t_{j,n}) \int_{t_{j,n}}^{t_{j+1,n}} \left( \int_{t_{j,n}}^{t} \sigma(s) \, dW_s \right) \sigma(t) \, dW_t =: A_1 + A_2,$$

where

$$A_1 := \sum_{j=0}^{k_n-1} (h_N(t_{j,n}) - h(t_{j,n})) \left[ 2 \int_{t_{j,n}}^{t_{j+1,n}} \left( \int_{t_{j,n}}^{t} \sigma(s) \, dW_s \right) \sigma(t) \, dW_t \right]$$

and

$$A_2 := \sum_{j=0}^{k_n-1} h(t_{j,n}) \left[ 2 \int_{t_{j,n}}^{t_{j+1,n}} \left( \int_{t_{j,n}}^{t} \sigma(s) \, dW_s \right) \sigma(t) \, dW_t \right].$$

Consider the term $A_2$. Using Theorem 5.5 in [27] and generalization in Proposition 1 in [37], we have that

$$(\rho(n))^{-1/2} \sum_{j=0}^{k_n-1} h(t_{j,n}) \left[ 2 \int_{t_{j,n}}^{t_{j+1,n}} \left( \int_{t_{j,n}}^{t} \sigma(s) \, dW_s \right) \sigma(t) \, dW_t \right]$$

converges to

$$\int_0^{2\pi} \sqrt{2H'(t)} h(t) \sigma^2(t) \, dB_t,$$

where $B(t)$ is a Brownian motion independent from $p(t)$ and the convergence is in law stably (for the definition of stable convergence see Section 2 in [27]). Therefore as $\rho(n) \to 0$ the limit of

$$(\rho(n))^{-1/2} \sum_{j=0}^{k_n-1} h(t_{j,n}) \left[ 2 \int_{t_{j,n}}^{t_{j+1,n}} \left( \int_{t_{j,n}}^{t} \sigma(s) \, dW_s \right) \sigma(t) \, dW_t \right]$$

is mixing Gaussian with a random variance equal to $2 \int_0^{2\pi} H'(t) h^2(t) \sigma^4(t) \, dt$ and the Gaussian variable is independent from $p(t)$. Finally consider the term $A_1$. By (13) and Burkholder–Davis–Gundy inequality we have

$$(\rho(n))^{-1/2} E\left[ \left( \sum_{j=0}^{k_n-1} (h_N(t_{j,n}) - h(t_{j,n})) \left( \int_{t_{j,n}}^{t_{j+1,n}} \left( \int_{t_{j,n}}^{t} \sigma(s) \, dW_s \right) \sigma(t) \, dW_t \right) \right)^2 \right]$$

$$\leq (\rho(n))^{-1/2} \left( \sup_{0 \leq t \leq 2\pi} |h_N(t) - h(t)| \right)^2$$

$$\times E\left[ \left( \sum_{j=0}^{k_n-1} \int_{t_{j,n}}^{t_{j+1,n}} \left( \int_{t_{j,n}}^{t} \sigma(s)\, dW_s \right) \sigma(t)\, dW_t \right)^2 \right]$$

$$\leq (\rho(n))^{-1} \left( \omega_h\left(\frac{4}{N}\right) \right)^2 CE\left[ \int_0^{2\pi} \sigma^4(t)\, dt \right],$$

where $C$ is a suitable constant. Therefore $(\rho(n))^{-1/2} A_1$ goes to 0 in probability as $n \to \infty$. $\square$

PROOF OF THEOREM 4.5. The proof is similar to the univariate case, except it is necessary to consider nonsynchronous observations for different assets. With the notation in Theorem 3.3 we have to consider

$$(40) \quad \sum_{|k|<N} \left( 1 - \frac{|k|}{N} \right) c_k(h) \overline{\alpha_k}(N, p_n^1, p_n^2) - \sum_{|k|<N} c_k(h) \overline{c_k}(\Sigma^{12}).$$

We prove that

$$(\rho(n))^{-1/2} \sum_{|k|<N} \left( 1 - \frac{|k|}{N} \right) c_k(h) \overline{\alpha_k}(N, p_n^1, p_n^2)$$

$$= (\rho(n))^{-1/2} \sum_{|k|<N} \left( 1 - \frac{|k|}{N} \right) c_k(h)$$

(41)
$$\times \left( \sum_{ij} e^{ikt_j^2} D_N(t_j^2 - t_i^1) \int_{t_i^1 \vee t_j^2}^{t_{i+1}^1 \wedge t_{j+1}^2} dp^1(t) \int_{t_i^1 \vee t_j^2}^{t_{i+1}^1 \wedge t_{j+1}^2} dp^2(t) \right)$$

$$+ o_p(1).$$

In fact, we split

$$\overline{\alpha_k}(N, p_n^1, p_n^2)$$

(42)
$$= \frac{1}{2\pi} \frac{1}{2N+1} \sum_{|s| \leq N} \sum_{i,j} e^{ikt_j^2} e^{-is(t_j^2 - t_i^1)} \delta_{I_i^1}(p^1) \delta_{J_j^2}(p^2) I_{\{I_i^1 \cap J_j^2 \neq \varnothing\}}$$

(43)
$$+ \frac{1}{2\pi} \frac{1}{2N+1} \sum_{|s| \leq N} \sum_{i,j} e^{ikt_j^2} e^{-is(t_j^2 - t_i^1)} \delta_{I_i^1}(p^1) \delta_{J_j^2}(p^2) I_{\{I_i^1 \cap J_j^2 = \varnothing\}}.$$

Consider (42). Observe that

$$(44) \quad I_{\{I_i^1 \cap J_j^2 \neq \varnothing\}} = I_{\{t_i^1 \vee t_j^2 < t_{i+1}^1 \wedge t_{j+1}^2\}}.$$

Therefore

$$\delta_{I_i^1}(p^1)\delta_{J_j^2}(p^2)I_{\{I_i^1 \cap J_j^2 \neq \varnothing\}}$$

$$(45) \quad = \Bigg\{ \int_{t_i^1 \vee t_j^2}^{t_{i+1}^1 \wedge t_{j+1}^2} dp^1(s) \int_{t_i^1 \vee t_j^2}^{t_{i+1}^1 \wedge t_{j+1}^2} dp^2(s)$$

$$(46) \quad + \left(\int_{t_i^1}^{t_{i+1}^1} dp^1(s)\right)\left(\int_{t_j^2}^{t_i^1 \vee t_j^2} dp^2(s) + \int_{t_{i+1}^1 \wedge t_{j+1}^2}^{t_{j+1}^2} dp^2(s)\right)$$

$$(47) \quad + \left(\int_{t_i^1 \vee t_j^2}^{t_{i+1}^1 \wedge t_{j+1}^2} dp^2(s)\right)\left(\int_{t_i^1}^{t_i^1 \vee t_j^2} dp^1(s) + \int_{t_{i+1}^1 \wedge t_{j+1}^2}^{t_{i+1}^1} dp^1(s)\right)\Bigg\}$$

$$\times I_{\{t_i^1 \vee t_j^2 < t_{i+1}^1 \wedge t_{j+1}^2\}}.$$

We first consider (46). The terms (43) and (47) are estimated analogously.

By splitting $I_{\{t_i^1 \vee t_j^2 < t_{i+1}^1 \wedge t_{j+1}^2\}}$ into the cases $I_{\{t_i^1 < t_j^2 < t_{j+1}^2 < t_{i+1}^1\}}$, $I_{\{t_i^1 < t_j^2 < t_{i+1}^1 < t_{j+1}^2\}}$, $I_{\{t_j^2 < t_i^1 < t_{i+1}^1 < t_{j+1}^2\}}$ and $I_{\{t_j^2 < t_i^1 < t_{j+1}^2 < t_{i+1}^1\}}$, we obtain:

$$\sum_{|k|<N}\left(1 - \frac{|k|}{N}\right)c_k(h)\sum_{i,j}e^{ikt_j^2}D_N(t_j^2 - t_i^1)\left(\int_{t_i^1}^{t_{i+1}^1} dp^1\right)$$

$$\times \left(\int_{t_j^2}^{t_i^1 \vee t_j^2} dp^2(s) + \int_{t_{i+1}^1 \wedge t_{j+1}^2}^{t_{j+1}^2} dp^2(s)\right)I_{\{t_i^1 \vee t_j^2 < t_{i+1}^1 \wedge t_{j+1}^2\}}$$

$$= 2\sum_{ij} h_N(t_j^2)D_N(t_j^2 - t_i^1)\left(\int_{t_i^1}^{t_{i+1}^1} dp^1(s)\right)\left(\int_{t_j^2}^{t_i^1} dp^2(s) + \int_{t_{i+1}^1}^{t_{j+1}^2} dp^2(s)\right)$$

$$\times I_{\{t_j^2 < t_i^1 < t_{i+1}^1 < t_{j+1}^2\}}.$$

Define

$$U(\phi_1, \phi_2) := \sum_{i,j} h_N(t_j^2)D_N(t_i^1 - t_j^2)I_{[t_i^1, t_{i+1}^1[}(\phi_1)I_{[t_j^2, t_i^1[}(\phi_2),$$

where $D_N(s)$ is the rescaled Dirichlet kernel. Then

$$\sum_{i,j} h_N(t_j^2)D_N(t_i^1 - t_j^2)\int_{t_i^1}^{t_{i+1}^1} dp^1(s)\int_{t_j^2}^{t_i^1} dp^2(s)I_{\{t_j^2 < t_i^1 < t_{i+1}^1 < t_{j+1}^2\}}$$

$$= \int\int_{\phi_2 < \phi_1} U(\phi_1, \phi_2)\,dp^1(\phi_1)\,dp^2(\phi_2).$$

By Itô energy identity and Plancherel equality, we have

$$E\left[\left(\int\int_{\phi_1 < \phi_2} U(\phi_1, \phi_2)\,dp^1(\phi_1)\,dp^2(\phi_2)\right)^2\right]$$

$$= E\left[\int_0^{2\pi}\left(\int_0^{\phi_2} U(\phi_1,\phi_2)\,dp^1(\phi_1)\right)^2 \Sigma^{22}(\phi_2)\,d\phi_2\right]$$

$$\leq C_1 E\left[\int_0^{2\pi}\int_0^{2\pi} D_N^2(\phi_1-\phi_2)h_N^2(\phi_1)d\phi_1 d\phi_2\right]$$

$$\times \operatorname{ess\,sup}\|\Sigma^{22}\|_{L^\infty}\operatorname{ess\,sup}\|\Sigma^{11}\|_{L^\infty}$$

$$\leq C_2 \frac{2\pi}{2N+1}\operatorname{ess\,sup}\|\Sigma^{22}\|_{L^\infty}\operatorname{ess\,sup}\|\Sigma^{11}\|_{L^\infty},$$

where $C_1$ and $C_2$ are suitable constants. It follows

$$\rho(n)^{-1/2} E\left[\left|\sum_{i,j} h_N(t_j^2) D_N(t_i^1-t_j^2)\int_{t_i^1}^{t_{i+1}^1} dp^1(s)\int_{t_j^2}^{t_i^1} dp^2(s)I_{\{t_j^2<t_i^1<t_{i+1}^1<t_{j+1}^2\}}\right|^2\right]$$

$$\leq C\rho(n)^{-1/2}N^{-1}.$$

Finally, we have to consider the term (45). Define

$$G_N := \sum_{ij} h_N(t_j^2)(D_N(t_j^2-t_i^1)-1)\int_{t_i^1\vee t_j^2}^{t_{i+1}^1\wedge t_{j+1}^2} dp^1(t)\int_{t_i^1\vee t_j^2}^{t_{i+1}^1\wedge t_{j+1}^2} dp^2(t).$$

We have

$$E\left[\left(\sum_{ij} h_N(t_j^2)(D_N(t_j^2-t_i^1)-1)\int_{t_i^1\vee t_j^2}^{t_{i+1}^1\wedge t_{j+1}^2} dp^1(t)\int_{t_i^1\vee t_j^2}^{t_{i+1}^1\wedge t_{j+1}^2} dp^2(t)\right)^2\right]$$

$$\leq \rho(n)^2 N^2 \operatorname{ess\,sup}\|\Sigma^{22}\|_{L^\infty}\operatorname{ess\,sup}\|\Sigma^{11}\|_{L^\infty}.$$

It follows

$$\rho(n)^{-1/2}E[|G_N|^2]\leq C\rho(n)^{3/2}N^2,$$

which goes to 0 by the hypothesis. Therefore (40) is equal to

$$\sum_{ij} h_N(t_j^2)\int_{t_i^1\vee t_j^2}^{t_{i+1}^1\wedge t_{j+1}^2} dp^1(t)\int_{t_i^1\vee t_j^2}^{t_{i+1}^1\wedge t_{j+1}^2} dp^2(t) - \int_0^{2\pi} h_N(t)\Sigma^{12}(t)\,dt + o_p(1).$$

Then the proof proceeds as in Theorem 4.2. $\square$

**Acknowledgments.** We heartily thank the Editors, an anonymous Associate Editor and two anonymous referees whose remarks led to a substantial improvement of this article.

ACADÉMIE DES SCIENCES
INSTITUT DE FRANCE
RUE SAINT LOUIS EN L'ISLE
75004 PARIS
FRANCE
E-MAIL: paul.malliavin@upmc.fr

DEPARTMENT OF MATHEMATICS FOR DECISIONS
UNIVERSITY OF FIRENZE
VIA LOMBROSO 6/17
FIRENZE
ITALY
E-MAIL: mariaelvira.mancino@dmd.unifi.it
URL: http://www.dmd.unifi.it